\documentclass[eprint,
]{actapoly}

\usepackage{amsthm}

\usepackage[utf8]{inputenc}

\begin{document}

\title[Some remarks on a mathematical model for water flow in porous media]
{Some remarks on a mathematical model for water flow in porous media with competition between transport and diffusion}

\author[J. Runcziková]{Judita Runcziková}{mat}
\author[J. Chleboun]{Jan Chleboun}{mat}
\correspondingauthor[C. Gavioli]{Chiara Gavioli}{mat}{chiara.gavioli@cvut.cz}
\author[P. Krejčí]{Pavel Krejčí}{mat}

\institution{mat}{Czech Technical University in Prague, Faculty of Civil Engineering, Department of Mathematics, Thákurova~7, 166~29 Prague, Czech Republic}

\begin{abstract}
The contribution deals with the mathematical modelling of fluid flow in porous media, in particular water flow in soils, with the aim of describing the competition between transport and diffusion. The analysis is based on a mathematical model developed by B.~Detmann, C.~Gavioli, and P.~Krejčí, in which the effects of gravity are included in a novel way. The model consists of a nonlinear partial differential equation describing both the diffusion and the gravitational transport of water. Although analytical solutions can be obtained for some special cases, only numerical solutions are available in more general situations. The solving algorithm is based on a time discretisation and the finite element method, and is written in Matlab. The results of the numerical simulations are shown and the behaviour of the model is discussed.

\end{abstract}

\keywords{porous media, water flow, transport, diffusion, PDE, numerical methods}

\maketitle

\section{Introduction}
Mathematical models for water flow in unsaturated porous media, such as soils, play an important role in many applications. They have proven to be useful in geoscience and environmental engineering (e.g.~spread of contamination into soil and groundwater, and their sanitation), in civil engineering (e.g.~the effect of water content on structures), and in other industrial applications.

A simple mathematical model for water flow in sa\-tu\-ra\-ted porous media is obtained from the mass balance equation (i.e.~the continuity equation) combined with the Darcy law for the water mass flux. If the porous medium is unsaturated, the water flux is described by a saturation-dependent variant of the Darcy law (also called the Darcy-Buckingham law), and the resulting equation is known as the Richards equation \cite{Cerny}. Many numerical methods have been used to solve the Richards equation. Some of them, with their advantages and limitations, are studied for e\-xam\-ple in \cite{Haverkamp, Ross, List, Farthing, Kamil}.
In \cite{Gandolfi} a comparative analysis is made between vertical water flow models derived from numerical solutions of the Richards equation and reservoir cascade schemes used in hydrological modelling. Additionally, dual-porosity models offer an alternative framework for simulating water flow in unsaturated soils, as exemplified by \cite{Lewandowska}.

While this paper focuses on a macroscopic-scale model examining the competition between transport and diffusion, it is worth noting that the competitive roles of driving forces in porous media have also been explored at the microscopic level, for instance for the pore-network model (discrete model) in \cite{Mohammadi}.

More elaborate mathematical models incorporating additional effects can be found in the literature. Pre\-fe\-ren\-tial flow (i.e.~pre\-fe\-ren\-tial paths through which the water flows more easily \cite{Allaire}), observed for e\-xam\-ple in \cite{Babel}, has been addressed in \cite{Vogel} and in the models summarized in \cite{Gerke}. Furthermore, in \cite{Schweizer} the Richards equation with hysteresis is introduced. The effects of capillary hysteresis on finger flow, particularly within non-equilibrium Richards equation mo\-dels, are investigated in \cite{Roche1}. Finally, a semi-continuum model capturing both diffusion-like and finger-like flow patterns is shown in \cite{Kmec}.

This paper is structured as follows. In Section~2, the partial differential equation (PDE) that models unsaturated fluid flow in a porous medium is presented, and the solving algorithm is described. In Section~3, the results of numerical simulations are shown and discussed. Examples with different initial and boun\-da\-ry conditions related to elementary hydrological processes in a soil column are considered in this paper.

\section{Materials and methods}

\subsection{The mathematical model \label{mmodel}}
A mathematical model has been proposed by Bettina Detmann, Chiara Gavioli, and Pavel Krejčí in \cite{BDCGPK} to describe fluid flow in an unsaturated porous medium under the influence of both diffusion and gravitational transport. Its main feature is a "stickiness" condition: if the saturation stays below a certain threshold $\bar{s}$, no transport takes place, and water flows only by diffusion. This behaviour is encoded in \eqref{mode01} below, where the notation $(\,\cdot\,)^+$ re\-pre\-sents the nonnegative part of a function.

The model is a good approximation of the real behaviour when $s$ stays away from the (expected) maximum value $1$. The full saturation case could be taken into account by including a term of the form $(1-s)^+$ in the transport term. This will be the subject of future research.

The mathematical model derived from the mass balance equation is the following:
\begin{align}
s_t - \kappa \Delta p - 2 \alpha g \left(s - \bar{s} \right)^+ s_z = 0, 
\label{mode01}
\end{align}
where $s = s({\bf x},t) \in (0,1)$ is the saturation (i.e.~the ratio between the water volume and pore volume), $s_t$ its time derivative, $s_z$ its derivative with respect to $z$
(i.e.~the variable representing the vertical direction), and $p$ is the capillary pressure. The time is denoted by $t \geq 0$, the porous body is represented by a domain $\Omega \subset \mathbb{R}^3$ with coordinates ${\bf x} = (x,y,z) \in \Omega$ and with Lipschitz boun\-da\-ry $\partial \Omega$. The parameter $\kappa$ is the coefficient, which corresponds to the permeability of the medium, $\alpha$ is the characteristic time related to the friction at the liquid-solid interface, $g$ is the gravity constant, and $\bar{s} \in (0,1)$ denotes the critical saturation value. If $s\le \bar{s}$, then $s$ can increase or decrease only due to the flow driven by diffusion. The model is complemented by a constitutive relation between the saturation $s$ and the capillary pressure $p$. Experimental evidence (see e.g.~\cite{Poulovassilis} for tests in a close-to-natural environment) indicates that the dependence is of hysteresis type. However, in the case of a purely hysteretic pressure-saturation relationship, the degeneracy of the problem could affect the convergence of the numerical scheme, see \cite{DCDS}. For simplicity, and as a first step towards more complex models, only linear dependence between $s$ and $p$ is considered in this paper.

The model in \cite{BDCGPK} is also coupled with an initial condition
\begin{align}
s({\bf x},0)=s_0({\bf x})
\end{align}
and boundary conditions assuming no flow through the vertical parts of the boundary
\begin{align}\label{bou1}
Q \cdot n = 0,
\end{align}
where $Q$ is the water mass flux vector, $n$ the unit outward normal vector, and $\cdot $ denotes the scalar pro\-duct. The model is completed by boundary conditions prescribing the inflow/outflow through the horizontal parts of the boundary
\begin{align}\label{bou0}
Q \cdot n = \beta (s-s_{out}),
\end{align}
where $\beta$ depends on the permeability coefficient and $s_{out}$ is the saturation outside the domain. 
In what follows, other boundary conditions on the horizontal part of the boundary will also be taken into account for this model, including the Dirichlet type for its numerical simplicity.

The model equation, the initial condition, and the boundary conditions for the 1-D flow are listed in Section~\ref{Algorithm}. More details can be found in \cite{BDCGPK}.

\subsection{Algorithm} \label{Algorithm}
The vertical flow of water through a column of soil is modelled as 1-D flow through an unsaturated porous body. Transport and diffusion effects are taken into account. The relation between pressure and saturation is chosen to be li\-near, that is, 
\begin{align} \label{linear}
s = \gamma p,
\end{align}
where $\gamma = 1$.

Equation \eqref{mode01} in the 1-D case coupled with relation (\ref{linear}) reduces to
\begin{align}
s_t - \kappa s_{zz} - 2 \alpha g \left(s - \bar{s} \right)^+  s_z  = 0,
\label{mode02}
\end{align}
where $\kappa >0$, $\alpha>0$, $g>0$, and $\bar{s} \in (0,1)$ are constants. The domain is defined as $\Omega = (-h,0)$ with coordinate $z \in (-h,0)$, so that $s = s(z,t)$ and $p = p(z,t)$. 

The equation-solving algorithm, based on the finite element method, is developed in Matlab. First, pa\-ra\-me\-ter values, time interval, and spatial mesh are set up, as well as the initial and boundary conditions. Then, the \texttt{pdepe} procedure (namely, the solver for systems of parabolic and elliptic PDEs of one spatial variable $x$ and time $t$) is used. The options of the solver are set to their default values, except the relative error tolerance (\texttt{'RelTol'}) which in the examples in Section~\ref{results} is set equal to $10^{-5}$. The domain is discretised by fixed division points, the spatial mesh is prescribed by matrices of $x$-coordinates and of $y$-coordinates and of the numbering of these points. The time mesh is set up automatically, because the time is discretised by the \texttt{ode15s} solver. This multistep solver is a variable-step, variable-order solver based on the numerical differentiation formulas of orders 1 to 5 \cite{matlab-ode}. More details about the \texttt{pdepe} solver can be found in \cite{procedura}. Then, the graphical outputs are generated.

The initial condition (i.e.~the saturation at the initial time $t=0$) is defined for all $z \in (-h,0)$ as
\begin{align}
s(z,0)=s_0(z),
\end{align}
where $s_0$ is a given function in $L^2(-h,0)$ taking values between 0 and 1.

The boundary conditions (BC) can be defined for $z=-h,0$ and $t >0$ as
\begin{itemize}
\item given values on the boundary (Dirichlet BC)
\begin{align} \label{bc01}
s(0,t)=f_+(t), \quad s(-h,t)=f_{-}(t);
\end{align}
\item given mass flux through the boundary (Neumann BC)
\begin{align} \label{bc02}
\kappa s_z + \alpha g \left[ (s - \bar{s} )^+ \right]^2 = \left\{ \begin{array}{lll}
f_+  \textrm{ for $z =0$,}\\
f_{-}  \textrm{ for $z =-h$;}
\end{array} \right.
\end{align}
\item mass flux proportional to inner/outer saturation difference (Newton BC)
\begin{align} \label{bc03}
\kappa s_z + \alpha g \left[ (s - \bar{s} )^+ \right]^2 = \left\{ \begin{array}{lll}
-\beta_+ (s-s_+) \textrm{ for $z=0$,}\\
\beta_- (s-s_-) \textrm{ for $z=-h$;}
\end{array} \right.
\end{align}
\end{itemize}
where $f_+, f_-$ are given time dependent functions, $\beta_+>0, \beta_->0$ are the permeability coefficients of the top and of the bottom, respectively, and $s_+,s_- \in [0,1]$ are the saturation values outside the domain in the positive and negative axis directions, respectively.
BC as in \eqref{bc01} with $f_+=f_-=0$ is considered in Section~\ref{Example3}, BC as in \eqref{bc02} with $f_+=f_-=0$ is considered in Section~\ref{Example1} and Section~\ref{Example2}.

The validity of the mass conservation is monitored by computing the definite integral of the solution $s$ over the domain by trapezoidal numerical integration. In the case of no inflow/outflow, the value of the integral has to remain constant over time. In the other cases, the value of the integral should comply with the difference between the inflow and outflow of water through the boundary. Thus, the water mass balance can be checked during the execution of the Matlab program.

\section{Results} \label{results}

An example of redistribution after infiltration is presented in Section~\ref{Example1}. The process of redistribution is described in \cite{Kutilek} and is also mentioned in \cite{Youngs}, where hysteresis effects are additionally included.
The initial condition corresponds to sa\-tu\-ra\-ted soil in the top part of the column, within a given depth below the terrain, with a steep transition (almost a jump) to dry soil.
In Section~\ref{Example2}, the process of wetting the soil column with a water source located at the bottom of the column is examined.
In Section~\ref{Example3}, the effect of $\kappa$ and $\bar{s}$ on a similar setting is studied.

\subsection{Example 1}\label{Example1}

\begin{figure}[htbp]
\centering
\includegraphics[width=\linewidth]{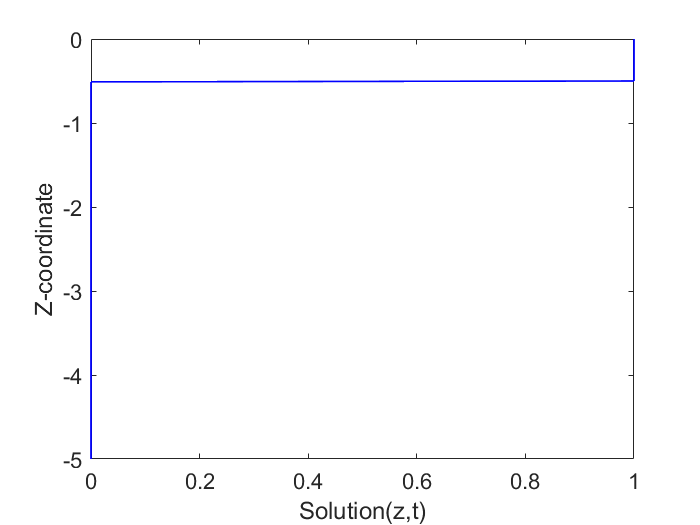}
\caption{Initial condition \eqref{ini1}.}
\label{fig:gr_PP2}
\end{figure}

\begin{figure}[htbp]
\centering
\includegraphics[width=\linewidth]{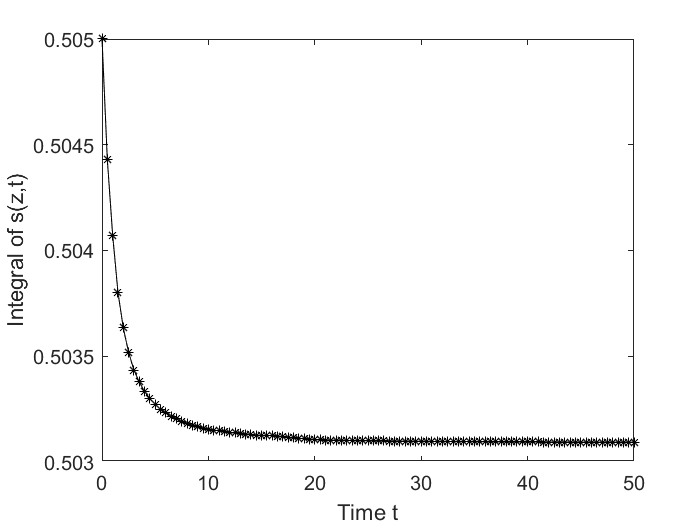}
\caption{Integral of the solution over the domain.}
\label{fig:gr_integral}
\end{figure}

In this first example, representing the redistribution after infiltration, the parameters are chosen as follows: permeability $\kappa = 0.005$, $2 \alpha g = 1$, $\bar{s}=0.25$, maximal depth $h = 5$, value of the space discretisation parameter $d = 0.01$. The value of the time step $\tau$ is chosen automatically by the solver. The initial condition (displayed in Figure~\ref{fig:gr_PP2}) is chosen as
\begin{align}\label{ini1}
s_0(z) = \left\{ \begin{array}{lll}
1 \textrm{ for $z>-0.50$,}\\
100z+51   \textrm{ for $z =[-0.51,-0.50]$,}\\
0 \textrm{ for $z<-0.51$,}
\end{array} \right.
\end{align}
which means saturated soil to the depth of 0.5 and dry soil below. The BC is that of impermeable boundary
\begin{align}\label{bou}
\kappa s_z + \alpha g \left[ (s - \bar{s} )^+ \right]^2 = 0,
\end{align}
so that the control integral of the saturation over the domain remains constant over time, as displayed in Figure~\ref{fig:gr_integral}. The exception is the small interval (of order $10^{-4}$) at the beginning, caused by a sharp change in the shape of the saturation curve.

\begin{figure}[htbp]
	\centering
	\includegraphics[width=7.8cm]{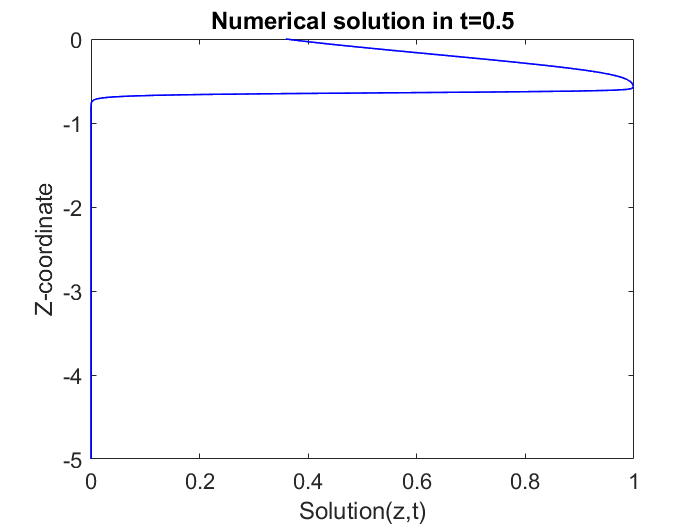}
	\\[3mm]
	\includegraphics[width=7.8cm]{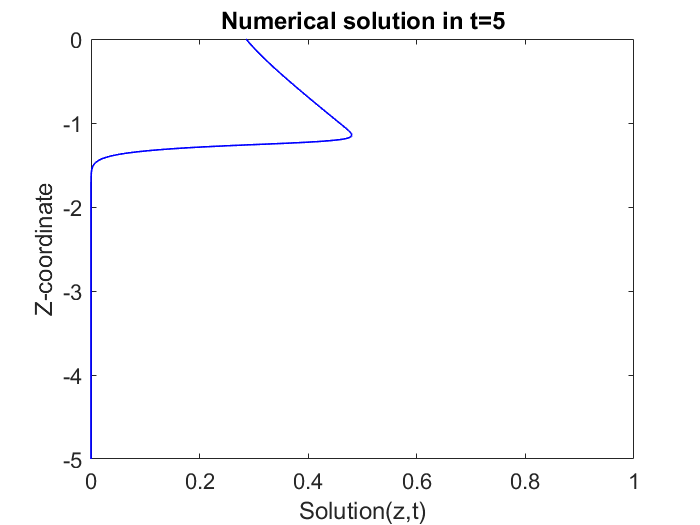}
	\\[3mm]
	\includegraphics[width=7.8cm]{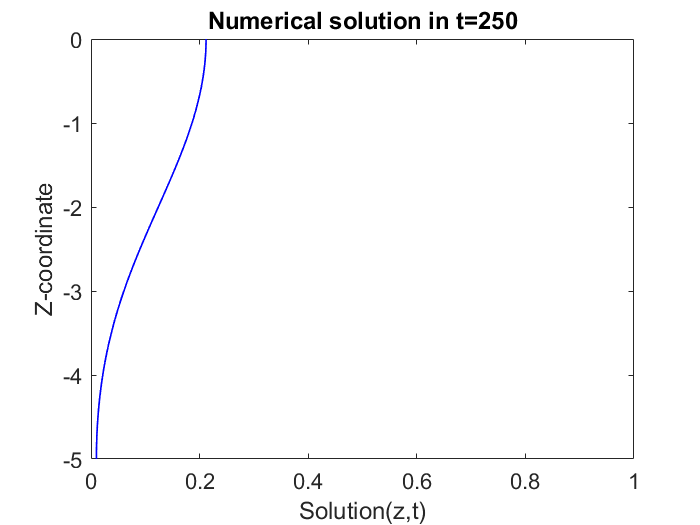}
	\\[3mm]
	\includegraphics[width=7.8cm]{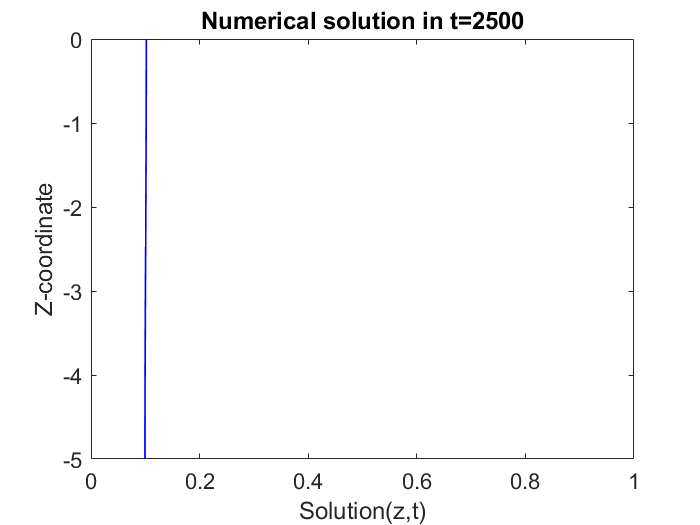}
	\caption{Evolution of the solution $s(z,t)$ over time with $\kappa=0.005$ and initial condition as in \eqref{ini1}.}
	\label{fig:gr_example}
\end{figure}

The saturation profile at selected times $t=0.5, 5, 250, 2500$ is plotted in Figure~\ref{fig:gr_example}.

At the beginning of the redistribution process, a rapid change in the saturation profile in the upper part can be noticed. A significant shift of the water content from the upper to the lower part is then observed, that is, the water content decreases in the upper part (above the wetting front) and increases in the lower part. The terrain saturation decreases from $s=1$ to near $s=\bar{s}$. This is caused by the dominant effect of transport in this zone. The effect of diffusion (note that $\kappa$ has a significantly large value) can be seen on the shape of the wetting front, which is slowly tilting and moving downward.

\begin{figure}[htbp]
\centering
\includegraphics[width=\linewidth]{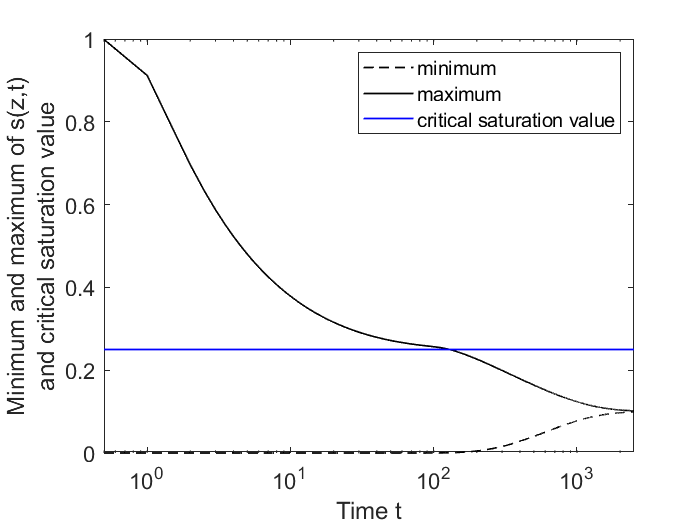}
\caption{Minimum and maximum of the solution $s(z,t)$ in the range defined by the time steps and the critical saturation value.}
\label{fig:gr_max}
\end{figure}

The saturation value over time decreases close to the terrain and increases in the lower part. According to \eqref{mode02}, if the value $s(z_i,t)$ of the saturation at a given depth $z_i$ is smaller than $\bar{s}$, then the transport term is not active at $z_i$, and the water movement is driven only by diffusion. Hence, if the saturation in the upper part decreases until it reaches the value $\bar{s}$, then only the diffusion term contributes to the redistribution of water. With the above parameter setting, the saturation along the entire soil column drops below $\bar{s}$ at time $t = 129$, as can be seen in Figure~\ref{fig:gr_max}. The water is then redistributed by diffusion until it reaches a uniform distribution profile along the soil column, and the maximal value of $s$ decreases until it meets the minimal value, as shown in Fi\-gu\-re~\ref{fig:gr_max}. From time $t = 609$, the difference between the maximum and minimum value of the solution is less than $0.1$. The process slows down as the saturation value decreases.

\subsection{Example 2}\label{Example2}

\begin{figure}[htbp]
	\centering
	\includegraphics[width=7.8cm]{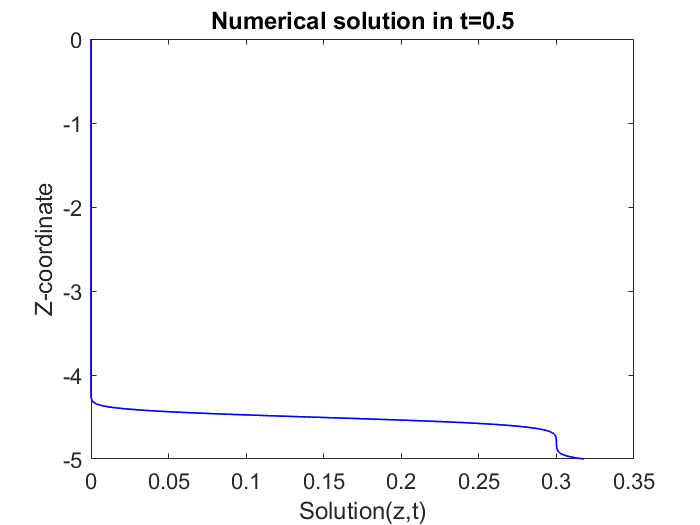}
	\\[3mm]
	\includegraphics[width=7.8cm]{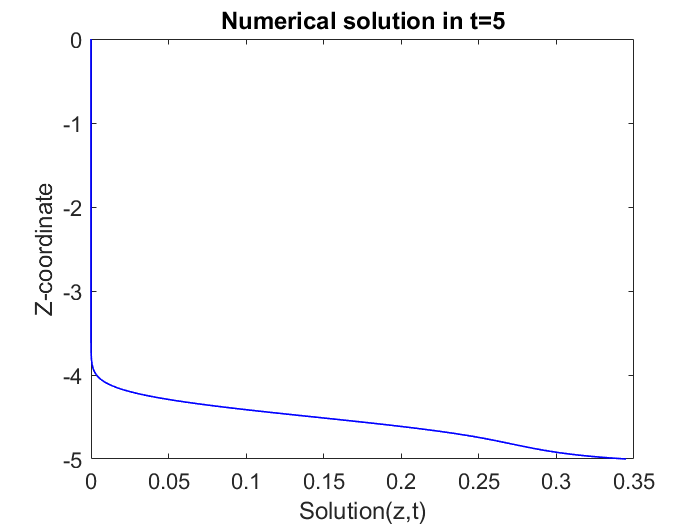}
	\\[3mm]
	\includegraphics[width=7.8cm]{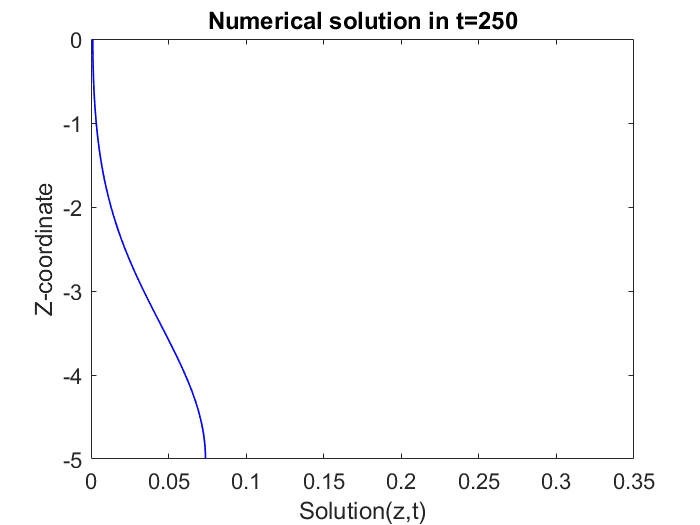}
	\\[3mm]
	\includegraphics[width=7.8cm]{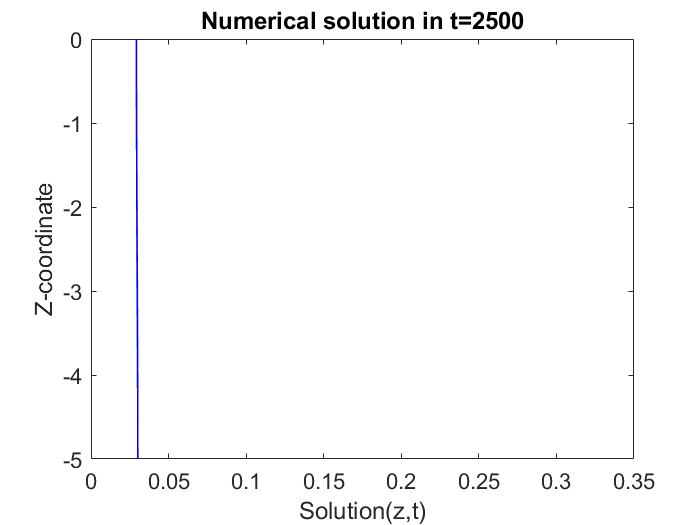}
	\caption{Evolution of the solution $s(z,t)$ over time  with $\kappa=0.005$ and initial condition as in \eqref{ini2}.}
	\label{fig:gr_example2}
\end{figure}

\begin{figure}[htbp]
\centering
\includegraphics[width=\linewidth]{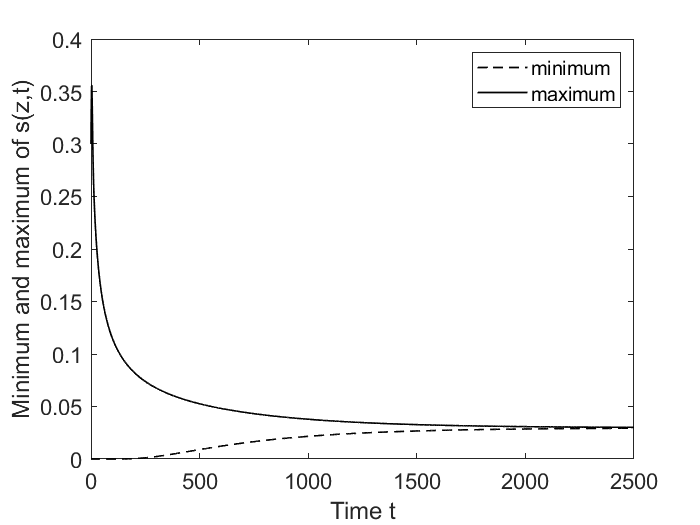}
\\[3mm]
\includegraphics[width=5.5cm]{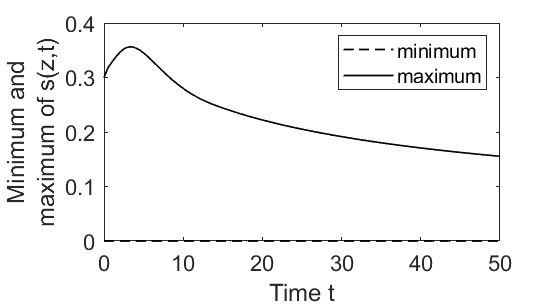}
\caption{Minimum and maximum of the solution $s(z,t)$ in the range defined by the time steps (top) and zoom to the beginning of the calculation (bottom).}
\label{fig:gr_max2}
\end{figure}

The initial condition is now, i.e.~in the example of soil wetting from the bottom water source, changed to
\begin{align}\label{ini2}
s_0(z) = \left\{ \begin{array}{lll}
0 \textrm{ for $z>-4.50$,}\\
-30z-135   \textrm{ for $z =[-4.50,-4.51]$,}\\
0.3 \textrm{ for $z<-4.51$,}
\end{array} \right.
\end{align}
which means unsaturated zone of thickness $0.49$ with saturation $s=0.3$ (hence $s>\bar{s}$) in the bottom part, with a steep transition to dry soil up to terrain. The other parameters and boundary conditions are as above. The saturation profile at selected times $t=0.5, 5, 250, 2500$ is plotted in Figure~\ref{fig:gr_example2}. First, flow driven by transport near the bottom and by diffusion above the unsaturated part is observed. Then the transport-driven flow in the lower part becomes weaker, the saturation starts to decrease and when the threshold value $\bar{s}$ is met, the transport term becomes inactive, so that only the upper part are continuously wetted by diffusion. The profiles of minimal and ma\-xi\-mal soil saturation over time are shown in Figure~\ref{fig:gr_max2}, where the peak in the bottom figure occurs at the time when the increase in saturation due to transport in the lower part ends.

\subsection{Example 3}\label{Example3}

\begin{figure}[htbp]
\centering
\includegraphics[width=\linewidth]{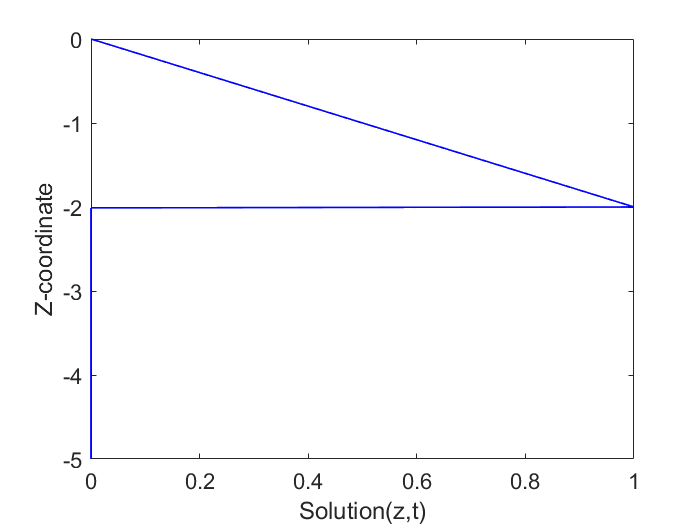}
\caption{Initial condition \eqref{ini3} with a steep wetting front.}
\label{fig:gr_PP}
\end{figure}

The influence of the choice of $\kappa$ and $\bar{s}$ on the sa\-tu\-ra\-tion profile is shown in the next examples. The following parameters are chosen: $\gamma=1$, $2 \alpha g = 1$, permeability $\kappa$ in the range $[0,0.01]$, critical saturation value $\bar{s}=0$, $\bar{s}=0.25$ and $\bar{s}=0.40$, maximal depth $h = 5$, value of the space discretisation parameter $d = 0.01$. The value of the time step $\tau$ is chosen automatically by the solver. The initial condition is now chosen as
\begin{align}\label{ini3}
s_0(z) = \left\{ \begin{array}{lll}
-0.5z \textrm{ for $z>-2.00$,}\\
100z+201   \textrm{ for $z =[-2.00,-2.01]$,}\\
0 \textrm{ for $z<-2.01$,}
\end{array} \right.
\end{align}
and is displayed in Figure~\ref{fig:gr_PP}.
For numerical simplicity, the boundary condition here is of Dirichlet type \eqref{bc01} with $f_+ = f_- = 0$.

The diffusion-driven flow is determined by the value of the permeability $\kappa$. If $\kappa=0$, the diffusion term in \eqref{mode02} is inactive. If the value of $\kappa$ is close to $0$, for example $\kappa=0.001$, the diffusion is negligible. Higher values of $\kappa$ are associated with strong diffusion.

\begin{figure}[htbp]
\centering
\includegraphics[width=\linewidth]{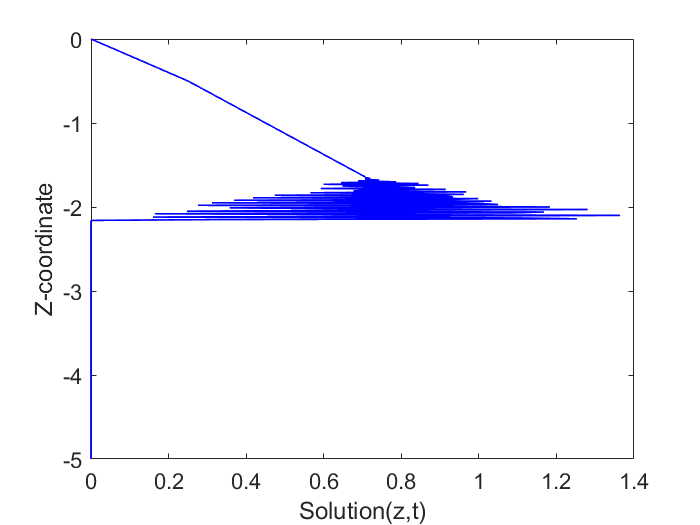}
\\[3mm]
\includegraphics[width=\linewidth]{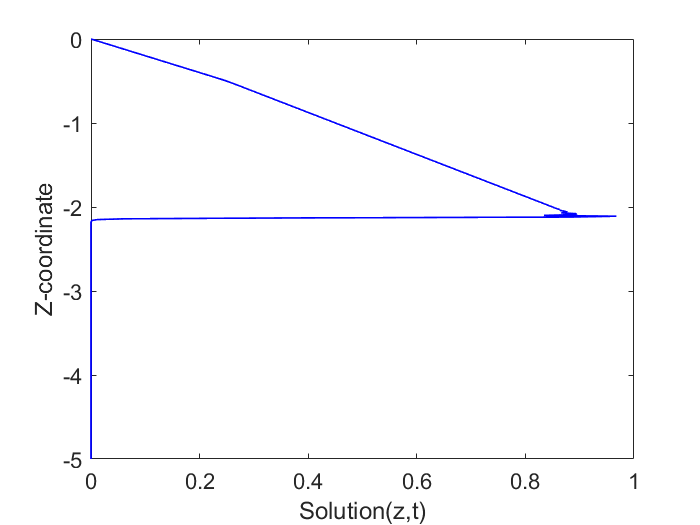}
\\[3mm]
\includegraphics[width=\linewidth]{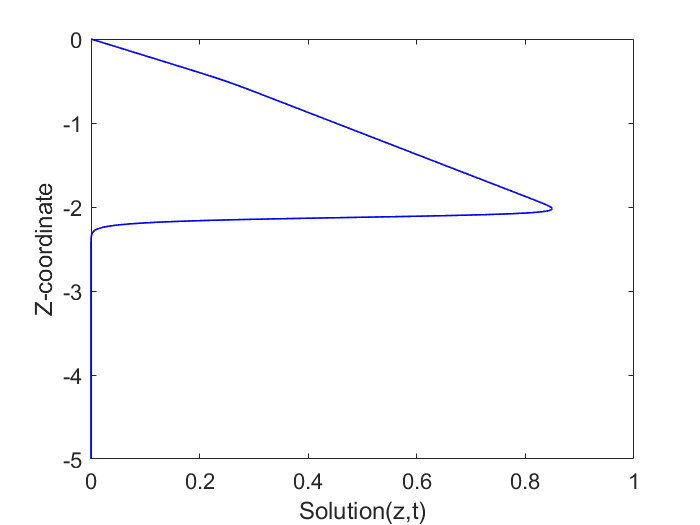}
\caption{Effect of the permeability $\kappa$ on the solution $s(z,0.5)$, with $\bar{s}=0.25$: without diffusion ($\kappa=0$, top), negligible diffusion ($\kappa=0.001$, centre), and strong diffusion ($\kappa=0.01$, bottom).}
\label{fig:gr_kappa}
\end{figure}

If the diffusion is strong enough, that is, the value of $\kappa$ is large enough, it is possible to choose an initial condition with a steep transition between two different saturation values, because the gradient of the solution soon becomes slightly smaller during the calculation due to the diffusion flow, and the simulation is stable. This situation occurs in this example in Fi\-gu\-re~\ref{fig:gr_kappa}, bottom, and also in Sections~\ref{Example1} and~\ref{Example2}.
Conversely, when the value of $\kappa$ is small, the gradient of the solution remains large, or even tends to infinity, and numerical instabilities (i.e.~oscillations of the solution) occur. Depending on the value of $\kappa$, these oscillations can cause calculation collapse, as shown in Figure~\ref{fig:gr_kappa}, top, or disappear, as in Fi\-gu\-re~\ref{fig:gr_kappa}, centre.

The same analysis applies if a high saturation gradient between the wet and dry parts of the porous body occurs later in the calculation.

\begin{figure}[htbp]
\centering
\includegraphics[width=\linewidth]{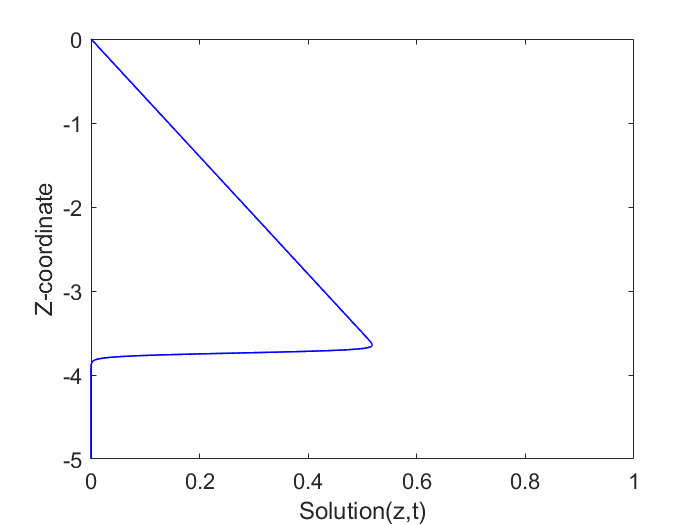}
\\[3mm]
\includegraphics[width=\linewidth]{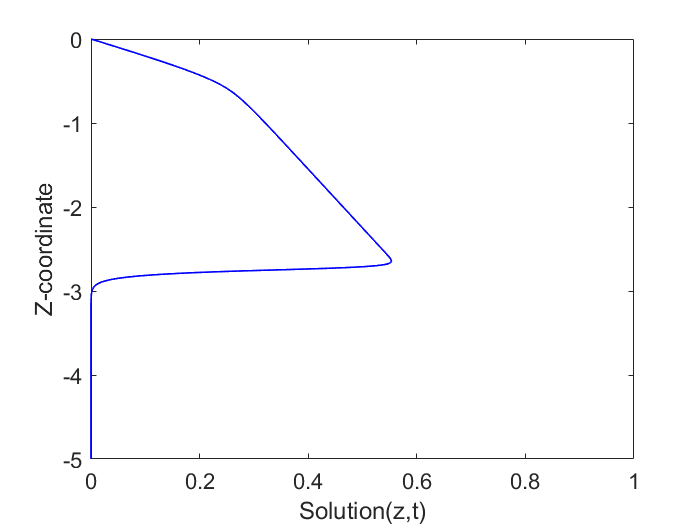}
\\[3mm]
\includegraphics[width=\linewidth]{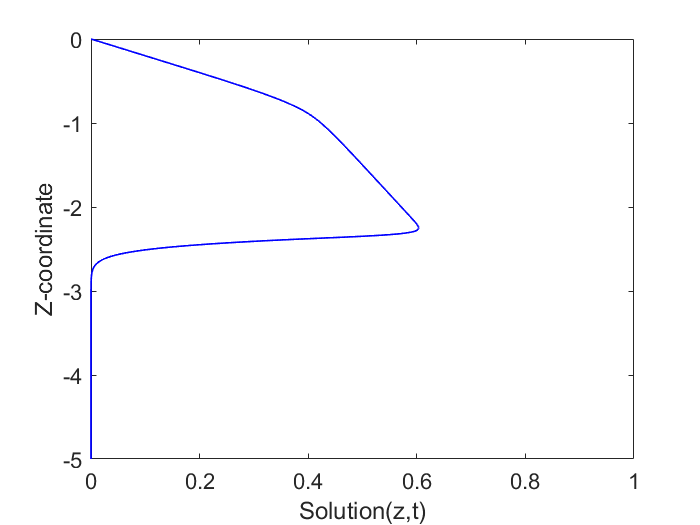}
\caption{Effect of $\bar{s}$ on the solution $s(z,5)$ with $\kappa=0.005$: without stickiness ($\bar{s} = 0$, top), with lower stickiness ($\bar{s}=0.25$, centre), and with higher stickiness ($\bar{s}=0.4$, bottom).}
\label{fig:gr_bars}
\end{figure}

In Figure~\ref{fig:gr_bars} the effect of the critical saturation value $\bar{s}$ on the solution $s$ is shown. When $\bar{s}$ is higher, the wetting front reaches a smaller depth.

\section{Conclusions}

A mathematical model for the unsaturated flow of water in a porous medium was studied. 
Numerical simulations of a vertical 1-D water flow in an un\-sa\-tu\-ra\-ted soil column were performed, and the effect of two parameters, the permeability $\kappa$ and the critical saturation value $\bar{s}$, was studied. Their magnitude has a significant impact on the solution and in some cases on the stability of the calculation.

The case of redistribution of water in soil after infiltration was shown, and the results were compared with theoretical physical expectations from the technical literature. The resulting saturation profiles are in agreement with \cite{Kutilek}. The case of wetting from below due to pressure difference was also con\-si\-de\-red. The model has shown a good approximation of the physical behaviour because, with the no-flux boun\-da\-ry condition \eqref{bou}, the saturation values stay under control even if we initially allow $s=1$ in some zone of given thickness as in \eqref{ini1}.

Initially, the model can be employed for simulating laboratory experiments, such as water flow through soil columns. For application to regional water flow models within catchment areas, more complex model formulations are required.
Furthermore, the in\-fluen\-ce of hysteresis, which can be incorporated using the Preisach model to describe the relationship between saturation and capillary pressure (similar to the approach in \cite{DCDS}), and preferential flow should be tho\-rough\-ly investigated.

From the modelling point of view, future research will include the case where the saturation $s$ approaches $1$. From the numerical point of view, investigations will be directed towards the de\-ve\-lop\-ment of an algorithm for the 2-D flow, to test the behaviour of the model in 2-D ($\Omega \subset \mathbb{R}^2$).

\begin{nomenclature}
\item[ ]{s}{Saturation}
\item[ ]{p}{Capillary pressure}
\item[ ]{\kappa}{Permeability}
\item[ ]{\bar{s}}{Critical saturation value}
\item[ ]{\alpha}{Characteristic time}
\item[ ]{g}{Gravity constant}
\end{nomenclature}

\begin{acknowledgements}
Supported by the project Centre of Advanced Applied Sciences (CAAS) No. CZ.02.1.01/0.0/0.0/16 019/0000778, which is co-financed by the European Union; by the Grant Agency of the Czech Technical University in Prague, grant No. SGS23/092/OHK1/2T/11; by the European Union's Horizon Europe research and innovation programme under the Marie Skłodowska-Curie grant agreement No 101102708; by the Czech Science Foundation (GA\v CR) project No. 24-10586S.
\end{acknowledgements}

\bibliographystyle{actapoly}
\bibliography{biblio}

\end{document}